\documentclass[12pt]{article} 
\usepackage{amssymb}
\usepackage{amsmath}
\usepackage{latexsym}

\newtheorem{theorem}{Theorem}[section]
\newtheorem{proposition}[theorem]{Proposition}
\newtheorem{corollary}[theorem]{Corollary}
\newtheorem{lemma}[theorem]{Lemma}

\newtheorem{remark}[theorem]{Remark}
\newtheorem{example}[theorem]{Example}
\newtheorem{definition}[theorem]{Definition}

\newcommand{\half}{\frac{1}{2}}

\newcommand{\vac}{\mathbf {1}}
\newcommand{\delz}{\partial_{z}}

\begin{document}

\title{A Generalized Vertex Operator Algebra for Heisenberg Intertwiners}
\author{
Michael P. Tuite
 and Alexander Zuevsky\thanks{
Supported by a Science Foundation Ireland Research Frontiers Grant, and
by Max-Planck-Institut f\"{u}r Mathematik, Bonn}
\\
School of Mathematics, Statistics and Applied Mathematics, \\
National University of Ireland Galway \\
University Road, Galway, Ireland}

\maketitle
\begin{abstract}
We consider the extension of the Heisenberg vertex operator algebra by all its irreducible modules. 
We give an elementary construction for the intertwining vertex operators and show that they satisfy a complex parametrized generalized vertex operator algebra. 
We illustrate some of our results with the example of integral lattice vertex operator superalgebras. 
\end{abstract}
\newpage
\section{Introduction}
The concept of a Vertex Operator Algebra (VOA) was introduced in \cite{B} and  \cite{FLM} 
and is essentially a rigorous algebraic approach to chiral conformal field theory in physics. 
This paper is devoted to one of the most basic examples, namely 
the rank $l$ Heisenberg VOA $M^{l}$ i.e. the chiral CFT consisting of $l$ free bosons.
An important application is the construction of a VOA $V$ containing a Heisenberg subVOA $M^{l}$ where we decompose 
$V$ into irreducible $M^{l}$-modules or extend $V$ by  $V$-modules or twisted $V$-modules related to the Heisenberg structure.
Thus in \cite{DLM} it is demonstrated how to extend a simple VOA $V$ by $g$-twisted $V$-modules for automorphisms $g$ generated by Heisenberg vectors \cite{Li1}. 
It is also shown in \cite{DLM} that intertwining vertex operators constructed on the larger space form a Generalized VOA \cite{DL}.  
In this paper we consider the extension of the Heisenberg VOA $M^{l}$ to the space, $\cal{M}$, given by the direct sum of all irreducible $M^{l}$-modules. 
We give an elementary construction for the intertwining vertex operators on $\cal{M}$ and show that these operators satisfy a complex parametrized Generalized VOA of a more general type than that defined in \cite{DL}.
We illustrate some of our results with the example of lattice VOSA $V_{L}$ for even or odd integral lattice $L$ and 
consider the extension of $V_{L}$ by twisted sectors for automorphisms generated by Heisenberg elements.   

We begin in Section~\ref{Intertwiner} with definitions and some properties of VOA modules and twisted modules. 
We consider the creative intertwiner vertex operators for a VOA $V$-module $W$ 
i.e. operators whose modes map $V$ onto $W$. 
In particular, we describe properties for these operators that are very similar to those 
for standard vertex operators.
In Section~\ref{Heis} we consider the rank $l$ Heisenberg VOA $M^{l}$ with irreducible module $M_{\alpha}$ for $\alpha\in \mathbb{C}^l$.
We give an elementary construction of the creative intertwiner (very similar in structure to vertex operators for a lattice VOA). 
In Section \ref{GVOA} we construct a $\mathbb{C}$-parametrized Generalized VOA on ${\cal M}=\oplus_{\alpha\in \mathbb{C}^l}  M_{\alpha}$, 
the direct sum of all irreducible $M^{l}$-modules. 
We also discuss skew-symmetry and show that there exists a unique invertible invariant symmetric bilinear form on ${\cal M}$. 
Finally in Section~\ref{LVOSA} we consider the example of the lattice VOSA $V_{L}$ for an even or odd integral lattice $L$. 
Using the Generalized VOA structure we construct the $g$-twisted $V_{L}$-module for an automorphism $g$ generated by a Heisenberg vector and show that this isomorphic to Li's construction \cite{Li1}. 
Finally, we conclude with a generalization of one of the main results of ref.~\cite{DLM} for $V_{L}$.     

\section{Creative Intertwiners for a Vertex Operator Algebra}
\label{Intertwiner}
We begin with a brief review of aspects of Vertex Operator Algebras and their modules, see refs.~\cite{B}, \cite{FLM},   \cite{FHL},  \cite{K}, \cite{LL}, \cite{MT} for more details. 
In particular, we are interested in describing properties of creative intertwiners which can be proved by a suitable modification of standard results in VOA theory.

We define the standard formal series 
\begin{eqnarray}
\delta\left(\frac{x}{y}\right)&=&\sum_{n\in \mathbb{Z}} x^{n}y^{-n},
\label{delta}\\
(x+y)^{\kappa}&=&\sum_{m\ge 0}\binom{\kappa}{m}x^{\kappa-m}y^m,
\label{xykappa}
\end{eqnarray} 
for any formal variables $x,y,\kappa$ where $\binom{\kappa}{m}=\frac{\kappa(\kappa-1)\ldots (\kappa-m+1)}{m!}$.  
\begin{definition}
 \label{VOA} 
A Vertex Operator Superalgebra (VOSA) is determined by a quadruple 
$(V,Y,\mathbf{1},\omega )$ 
as follows: $V$ is a superspace $V=V_{\bar{0}}\oplus V_{\bar{1}}$
with parity $p(u)=0$ or $1$ for $u\in V_{\bar{0}}$ or $V_{\bar{1}}$ respectively. $V$ also has a $\half\mathbb{Z}$-grading with $V=\bigoplus _{r \in \half\mathbb{Z}}V_{r}$  
with $\dim V_{r}<\infty$ and $V_{r}=0$ for $r\ll 0$.
   $\vac\in V_{0}$ is the vacuum vector and 
$\omega\in V_{2}$ is the conformal vector with properties described below.
 
$Y$ is a linear map $Y:V\rightarrow \mathrm{End}(V)[[z,z^{-1}]]$ for formal
variable $z$ so that for any vector $u\in V$ we have a vertex operator  
\begin{equation}
Y(u,z)=\sum_{n\in \mathbb{Z}}u(n)z^{-n-1}.  
\label{Ydefn}
\end{equation}
The linear operators (modes) $u(n):V\rightarrow V$ satisfy creativity 
\begin{equation}
Y(u,z)\vac = u +O(z)
\label{create}
\end{equation}
and lower truncation 
\begin{equation}
u(n)v=0,
\label{lowertrun}
\end{equation}
for each $u,v\in V$ and $ n\gg 0$.
For the conformal vector $\omega$ 
\begin{equation}
Y(\omega ,z)=\sum_{n\in \mathbb{Z}}L(n)z^{-n-2},  \label{Yomega}
\end{equation}
where $L(n)$ satisfy the Virasoro relation for some central charge $c$ 
\begin{equation}
[\,  L(m),L(n)\, ]=(m-n)L(m+n)+\frac{c}{12}(m^{3}-m)\delta_{m,-n}.
\label{Virasoro}
\end{equation}
Each vertex operator satisfies the translation property 
\begin{equation}
Y(L(-1)u,z)=\delz Y(u,z).  
\label{YL(-1)}
\end{equation}
The Virasoro operator $L(0)$ provides the $\half\mathbb{Z}$-grading with $L(0)u=ru$ for 
$u\in V_{r}$ and with $r\in \mathbb{Z}+\half p(u)$. 
Finally, the vertex operators satisfy the Jacobi identity
\begin{eqnarray}
\notag
&& z_0^{-1}\delta\left( \frac{z_1 - z_2}{z_0}\right) Y (u, z_1 )Y(v , z_2)    
  -(-1)^{p(u,v)}z_0^{-1} \delta\left( \frac{z_2 - z_1}{-z_0}\right) Y(v, z_2) Y(u , z_1 ) 
\\
\label{VOAJac}
&& 
= z_2^{-1}    
\delta\left( \frac{z_1 - z_0}{z_2}\right)
Y \left( Y(u, z_0)v, z_2\right),  
\end{eqnarray} 
where $p(u,v)=p(u)p(v)$.
\end{definition}
\begin{remark}
$(V,Y,\mathbf{1},\omega )$ is called a Vertex Operator Algebra (VOA) when $V_{\bar{1}}=0$.
\end{remark}
\medskip

Amongst other properties, these axioms imply 
locality, associativity, commutativity and skew-symmetry:
\begin{eqnarray}
(z_{1}-z_{2})^m
Y(u,z_{1})Y(v,z_{2}) 
&=& (-1)^{p(u,v)}(z_{1}-z_{2})^m
Y(v,z_{2})Y(u,z_{1}),
\notag
\\
&&
\label{Local}\\
(z_{0}+z_{2})^n Y(u,z_{0}+z_{2})Y(v,z_{2})w &=& (z_{0}+z_{2})^n Y(Y(u,z_{0})v,z_{2})w,
\label{Assoc}\\
u(k)Y(v,z)-(-1)^{p(u,v)} Y(v,z)u(k)
&=& \sum_{j\ge 0}\binom{k}{j}
Y(u(j)v,z)z^{k-j},
\label{Comm}\\
Y(u,z)v &=& (-1)^{p(u,v)} e^{zL(-1)}Y(v,-z)u,
\label{skew}
\end{eqnarray}
for $u,v,w\in V$ and integers $m,n\gg 0$  \cite{FHL},  \cite{K},  \cite{MT}. 
\medskip

We define the notion of a $V$-module  \cite{FHL}, \cite{LL}, \cite{X}. 
\begin{definition}
\label{module}
A $V$-module for a VOSA $V$ is a pair $(W, Y_W)$ where $W$ is a $\mathbb{C}$-graded vector space
$W=\bigoplus\limits_{r\in \mathbb{C}}W_{r}$ with $\dim W_{r}<\infty$ and 
where $W_{r+n}=0$ for all $r$ and $n\ll 0$. 
$Y_{W}$ is a linear map $Y_{W}:V\rightarrow \mathrm{End}(W)[[z,z^{-1}]]$ defining a module vertex operator
\begin{equation}
Y_{W}(u,z)=\sum_{n\in \mathbb{Z}}u_{W}(n)z^{-n-1},  
\label{YWdefn}
\end{equation}
 for each $u\in V$ with modes $u_{W}:W\rightarrow W$.  
For the vacuum vector $Y_W(\vac, z) = \mathrm{Id}_W$ and for the conformal vector
\begin{equation}
Y_{W}(\omega ,z)=\sum_{n\in \mathbb{Z}}L_{W}(n)z^{-n-2},  \label{YWomega}
\end{equation} 
where $L_{W}(0)w=r w$ for $w\in W_{r}$. The module vertex operators satisfy the Jacobi identity:
\begin{eqnarray}
\notag
&& z_0^{-1}\delta\left( \frac{z_1 - z_2}{z_0}\right) Y_{W} (u, z_1 )Y_{W}(v , z_2)    
  -(-1)^{p(u,v)} \delta\left( \frac{z_2 - z_1}{-z_0}\right) Y_{W}(v, z_2) Y_{W}(u , z_1 ) 
\\
&& 
= z_2^{-1}    
\delta\left( \frac{z_1 - z_0}{z_2}\right)
Y_{W} \left( Y(u, z_0)v, z_2\right). \label{VOAJacW} 
\end{eqnarray}
\end{definition}
A $V$-module $W$ is irreducible if no proper, nonzero subspace is invariant under all $u_{W}(n)$.  
The above axioms imply that $L_{W}(n)$ of \eqref{YWomega} satisfies the Virasoro algebra \eqref{Virasoro} for the same central charge $c$ and that the translation property holds: 
\begin{eqnarray}
Y_{W}(L(-1)u,z)&=&\delz Y_{W}(u,z).
\label{YL(-1)W}
\end{eqnarray}

\medskip
We next define the notion of a twisted $V$-module. 
Let $g$ be a $V$-automorphism $g$ i.e. a linear map preserving $\vac$ and $\omega$ such that 
\begin{equation*}
g Y(v,z) g^{-1}=  Y(g v,z),
\end{equation*}
for all $v \in V$. We assume that $V$ can be decomposed into $g$-eigenspaces
\begin{equation*}
V= \oplus_{\rho\in \mathbb{C}} V^{\rho}, 
\end{equation*}
where $V^{\rho}$ denotes the eigenspace of $g$ with eigenvalue $e^{2\pi i \rho}$.
\begin{definition}\label{Vgmodule}
A $g$-twisted $V$-module for a VOSA $V$ is a pair $(W^{g}, Y_{g})$ 
where $W^{g}$ is a $\mathbb{C}$-graded vector space
$W^{g}=\bigoplus\limits_{r\in \mathbb{C}}W^{g}_{r}$ with $\dim W_{r}<\infty$ and where    $W_{r+n}=0$ for all $r$ and  $n\ll 0$. 
$Y_{g}$ is a linear map $Y_{g}:V\rightarrow \mathrm{End\ }W^{g}\{z\}$, the vector space of 
$\mathrm{End\ }W^{g}$-valued formal series in $z$ with arbitrary complex powers of $z$. Then for $v\in V^{\rho}$
\begin{equation*}
Y_{g}(v,z)= \sum_{n \in \rho+ \mathbb {Z}} v_{g}(n)z^{-n-1}, 
\end{equation*}
with $v_{g}(\rho+l)  w=0$ for $w\in W^{g}$ and $l \in \mathbb {Z}$ sufficiently large. 
For the vacuum vector $Y_{g}(\vac, z) = \mathrm{Id}_{W^{g}}$ and for the conformal vector
\begin{equation}
Y_{g}(\omega ,z)=\sum_{n\in \mathbb{Z}}L_{g}(n)z^{-n-2},  \label{Ygomega}
\end{equation} 
where $L_{g}(0)w=r w$ for $w\in W^{g}$. The $g$-twisted vertex operators satisfy the twisted Jacobi identity: 
\begin{eqnarray}
\nonumber
& &  z_0^{-1} \delta\left( \frac{z_1 - z_2}{z_0}\right) Y_{g}(u, z_1) Y_{g}(v, z_2)  
- (-1)^{p(u,v)} z_0^{-1} \delta\left( \frac{z_2 - z_1}{-z_0}\right) Y_{g}(v, z_2) Y_{g}(u, z_1)  
\\
\label{gJac}
&& \qquad =
 z_2^{-1} \left( \frac{z_1 - z_0}{-z_2}\right)^{-\rho} 
\delta\left( \frac{z_1 - z_0}{-z_2}\right) Y_{g}(Y(u, z_0)v, z_2), 
\end{eqnarray}
for $u \in V^{\rho}$.  
\end{definition}
This definition is an extension of the standard one for $g$ of finite order where $\rho\in\mathbb{Q}$ \cite{D}, \cite{DLM} or for $g$ unitary where $\rho\in \mathbb{R}$ \cite{DLinM}. These axioms  imply that $L_{g}(n)$ of \eqref{Ygomega} satisfies the Virasoro algebra \eqref{Virasoro} for the same central charge $c$ and that the translation property holds: 
\begin{eqnarray}
Y_{g}(L(-1)u,z)&=&\delz Y_{g}(u,z).
\label{YL(-1)g}
\end{eqnarray} 
\medskip

We lastly restrict ourselves to a VOA $V$ with a $V$-module $(W, Y_W)$ and introduce the notion of creative intertwining vertex operators that satisfy an intertwining algebra of type $\binom{W}{W\ V}$ in the terminology of ref.  \cite{FHL}.
\begin{definition}\label{inter}
A creative intertwining vertex operator ${\cal Y}$ for a VOA $V$-module $(W, Y_W)$ is defined by a linear map ${\cal Y}:W\rightarrow \mathrm{Lin}(V,W)[[z,z^{-1}]]$ with 
\begin{equation}
{\cal Y}(w,z)=\sum_{n\in \mathbb{Z}}w(n)z^{-n-1},  
\label{calYdefn}
\end{equation}
for each $w\in W$ with modes $w(n):V\rightarrow W$. The intertwining vertex operator satisfies creativity 
\begin{equation}
{\cal Y}(w,z)\vac = w +O(z),
\label{createW}
\end{equation}
for each $w\in W$
and lower truncation 
\begin{equation}
w(n)v=0,
\label{lowertrunW}
\end{equation}
 for each $v\in V$, $w\in W$ and $n\gg 0$. The intertwining vertex operators satisfy the Jacobi identity:
\begin{eqnarray}
\notag
&& z_0^{-1}\delta\left( \frac{z_1 - z_2}{z_0}\right) Y_{W} (u, z_1 ){\cal Y}(w , z_2)    
  -z_0^{-1} \delta\left( \frac{z_2 - z_1}{-z_0}\right) {\cal Y}(w, z_2) Y(u , z_1 ) 
\\
&& 
= z_2^{-1}    
\delta\left( \frac{z_1 - z_0}{z_2}\right)
{\cal Y}\left( Y_{W}(u, z_0)w, z_2\right), \label{VOAJacInt} 
\end{eqnarray} 
for all $u\in V$ and $w\in W$.
\end{definition}
\medskip

These axioms imply that the intertwining and module vertex operators satisfy the following form of translation, locality, associativity, commutativity and skew-symmetry: 
\begin{eqnarray}
{\cal Y}(L_{W}(-1)w,z)&=&\delz {\cal Y}(w,z),
\label{YL(-1)YInt}
\\
(z_{1}-z_{2})^m
Y_{W}(u,z_{1}){\cal Y}(w,z_{2}) 
&=& (z_{1}-z_{2})^m
{\cal Y}(w,z_{2})Y(u,z_{1}),
\label{LocalW}\\
(z_{0}+z_{2})^n Y_{W}(u,z_{0}+z_{2}){\cal Y}(w,z_{2})v &=& (z_{0}+z_{2})^n {\cal Y}(Y_{W}(u,z_{0})w,z_{2})v,
\label{AssocW}\\
u_{W}(k){\cal Y}(w,z) - {\cal Y}(w,z) u(k)
&=& \sum_{j\ge 0}\binom{k}{j}
{\cal Y}(u_{W}(j)w,z)z^{k-j},
\label{CommW}\\
{\cal Y}(w,z)v &=& e^{zL_{W}(-1)}Y_{W}(v,-z)w,
\label{skewW}
\end{eqnarray}
 for $u,v\in V$, $w\in W$ and integers $m,n\gg 0$.

\medskip
We may obtain creative intertwining versions of other standard theorems for VOAs. 
The proofs are omitted since they are easily modified versions of standard VOA methods e.g.  \cite{LL}, \cite{K}. 
Define an intertwining normal ordering for $u\in V$ and $w\in W$ 
\begin{equation}
:Y_{W}(u,z){\cal Y}(w,z): = Y_{W}(u,z)_{-}{\cal Y}(w,z)+{\cal Y}(w,z)Y(u,z)_{+},
\label{norder}
\end{equation}
where 
\begin{eqnarray}
&Y(u,z)_{-}=\sum\limits_{n\ge 0}u(n)z^{-n-1},
\ &Y(u,z)_{+}=\sum\limits_{n<0}u(n)z^{-n-1},
\label{Ypm}
\\
&Y_{W}(u,z)_{-}=\sum\limits_{n\ge 0}u_{W}(n)z^{-n-1},
\ &Y_{W}(u,z)_{+}=\sum\limits_{n<0}u_{W}(n)z^{-n-1}.
\label{YWpm}
\end{eqnarray}
Extracting the coefficient of $z_{0}^{-2}z_{1}^{-1}$ in the creative intertwiner Jacobi identity \eqref{VOAJacInt} we obtain 
\begin{proposition} \label{prop_no}
For $u\in V$ and $w\in W$ we have 
\begin{equation}
:Y_{W}(u,z){\cal Y}(w,z):={\cal Y}(u_{W}(-1)w,z).
\label{No2}
\end{equation}
\end{proposition}   
Thus $:Y_{W}(u,z){\cal Y}(w,z):$ is local in the sense of \eqref{LocalW}. We also obtain
\begin{theorem}[Goddard's Uniqueness Theorem for Intertwiners]
\label{Goddard}
Let ${\cal W}(z)\in \mathrm{Lin}(V,W)[[z,z^{-1}]]$ be local in the sense of \eqref{LocalW} i.e. for each $u\in V$
\begin{equation}
(z_{1}-z_{2})^m
Y_{W}(u,z_{1}){\cal W}(z_{2})= (z_{1}-z_{2})^m
{\cal W}(z_{2})Y(u,z_{1}),\label{Wlocal}
\end{equation}
for $m\gg 0$. Suppose that for some $w\in W$
\begin{equation}
{\cal W}(z)\vac = e^{zL_{W}(-1)}w,
\label{Wcreate}
\end{equation}
then ${\cal W}(z)={\cal Y}(w,z)$, the creative intertwiner.
\end{theorem}

\section{Heisenberg Intertwiners}
\label{Heis}
In this Section we explicitly construct the creative intertwining operators for the irreducible modules of the rank~$l$ Heisenberg VOA $M^{l}$ generated by $l$ weight one  Heisenberg vectors $a^{i},\ i=1,\ldots,l$ with modes obeying 
\begin{equation}
\label{heisnorm}
[\, a^{i}(n), a^{j}(m)\, ]= n \delta_{n, -m}\delta_{i,j}.  
\end{equation}
$M^{l}$ is spanned by the Fock basis
\begin{equation}
a^{i_{1}}(-k_{1})a^{i_{2}}(-k_{2})\ldots a^{i_{r}}(-k_{r})\vac,\ k_{i}> 0,
\label{Fock}
\end{equation}
for  Virasoro vector $\omega=\frac{1}{2}\sum_{i=1}^{l}a^{i}(-1)^2\vac$ with central charge $l$.
 
The irreducible modules for $M^{l}$, denoted by $M_{\alpha}=M^{l}\otimes e^{\alpha}$ (with $M^{l}\cong M_{0}$) are indexed by a complex $l$-tuple $\alpha=\{\alpha^{1},\ldots,\alpha^{l}\} \in \mathbb{C}^{l}$ with
\begin{eqnarray}
a^{i}_{M_{\alpha}}(0)(u\otimes e^{\alpha}) &=& \alpha^{i} (u \otimes e^{\alpha}),\label{a0alpha}\\
a^{i}_{M_{\alpha}}(n)(u\otimes e^{\alpha}) &=& (a^{i}(n)u) \otimes e^{\alpha}, \ n\neq 0, \label{analpha}
\end{eqnarray}
for $Y_{M_{\alpha}}(a^{i},z)=\sum_{n\in \mathbb{Z}}a^{i}_{M_{\alpha}}(n)z^{-n-1}$ and $u\in M^{l}$. From now on we will employ the standard abbreviations of writing $Y(u,z)$ in place of $Y_{M_{\alpha}}(u,z)$, $u(n)$ in place of $u_{M_{\alpha}}(n)$ and $u$ in place of $u \otimes e^{0}$.  

We next construct the creative intertwiner for ${\cal Y}(u\otimes e^{\alpha},z)$ for $u\otimes e^{\alpha}\in M_{\alpha}$ for all $\alpha\in \mathbb{C}^{l}$. Much of the discussion is similar to the standard construction of lattice vertex operators e.g.  \cite{FLM},  \cite{K}. We first introduce the standard operators $q^{i}$ conjugate to $a^{i}(0)$ 
\begin{equation}
\label{pcom1}
\left[\,  a^{i}(n), q^{j}  \,\right] =  \delta_{n,0}\delta_{i,j}, 
\end{equation}
and identify 
\begin{equation*}
u \otimes e^{\alpha}=e^{\alpha \cdot q}(u \otimes e^{0})=e^{\alpha\cdot q}u,
\end{equation*} 
where $\alpha\cdot q=\sum_{i=1}^{l}\alpha^{i}q^{i}$ for $\alpha\in \mathbb{C}^{l}$.
We also define 
\begin{eqnarray}
\label{Yalphapm}
Y_{\pm}(\alpha, z) &=& \exp \left(\mp \;  \sum_{n>0}\frac{\alpha(\pm \; n)}{n} z^{\mp n}\right), 
\end{eqnarray} 
where $\alpha(n)=\alpha\cdot a(n)$.
These operators obey the following  \cite{Li1},  \cite{DLM}
\begin{proposition} 
\label{conjprop1}
For all $\alpha\in \mathbb{C}^{l}$ and $u\in M^{l}$ we find
\begin{eqnarray}
Y_{+}(\alpha,z_1)Y_{-}(\beta,z_2) &=& \left(1-\frac{z_2}{z_1}  \right)^{\alpha\cdot \beta}Y_{-}(\beta,z_2) Y_{+}(\alpha,z_1),
\label{Ypmcom}
\\
Y(Y_{+}(\alpha,-z_1)u,z_1)Y_{-}(\alpha,z_2)&=& Y_{-}(\alpha,z_2)Y(Y_{+}(\alpha,-z_{1}+z_{2})u,z_{1}),
\label{Yuconjminus}
\\
Y_{+}(\alpha,z_{1})Y(u,z_{2})&=&Y(Y_{+}(\alpha,z_{1}-z_{2})u,z_{2})
Y_{+}(\alpha,z_{1}),\label{Yuconjplus}
\\
Y(Y_{-}(\alpha, z_1) u, z_2)Y_{+}(\alpha, z_2 ) &=& z_2^{-\alpha(0)} (z_{2}+z_{1}  )^{\alpha(0)} Y_{-}(-\alpha, z_2)
Y_{-}(\alpha, z_{1}+z_{2}) 
\notag
\\
 & & .Y(u, z_2)Y_{+}(\alpha, z_{2} + z_{1}). \label{YYconj}
\end{eqnarray}
\end{proposition}
We also have 
\begin{proposition} 
\label{conjprop2}
For all $\alpha\in \mathbb{C}^{l}$ and $u\in M^{l}$ we find
\begin{eqnarray}
e^{-\alpha \cdot q} e^{zL(-1)}e^{\alpha \cdot q}e^{-zL(-1)} &=& Y_{-}(\alpha, z),
\label{L1conj}
\\
e^{-\alpha \cdot q} e^{zL(1)}e^{\alpha \cdot q}e^{-zL(1)} &=& Y_{+}(\alpha, \frac{1	}{z}),
\label{Lplus1conj}
\\
e^{-\alpha \cdot q} Y(u,z)e^{\alpha \cdot q} &=& Y(Y_{+}(\alpha, -z)u,z).
\label{Yuconj}
\end{eqnarray}
\end{proposition}
\noindent \textbf{Proof.}
From \eqref{pcom1} it follows that $[\, L(-1),q^{i}\, ]=a^{i}(-1)$. Hence we find
\begin{eqnarray*}
e^{-\alpha \cdot q} e^{zL(-1)}e^{\alpha \cdot q} &=& e^{z(L(-1)+\alpha(-1))}
\\
&=& Y_{-}(\alpha,z) e^{zL(-1)},
\end{eqnarray*}
from refs.  \cite{Li1},  \cite{DLM} giving \eqref{L1conj}. \eqref{Lplus1conj} follows similarly.

To prove \eqref{Yuconj} we first show that 
\begin{equation}
\left [\, q^{i}, Y(u,z)\, \right] = Y\left(X^{i}_{+}(z)u,z\right),
\label{Yuconj2}
\end{equation}
where $X^{i}_{+}(z)=\sum_{n>0} \frac{a^{i}(n)}{n} (-z)^n$. Assume that \eqref{Yuconj2} holds for every Fock vector $v$ with $m$ Heisenberg modes and consider $u=a^{j}(-k-1)v$ for $k\ge 0$. Then \eqref{pcom1} gives
\begin{eqnarray*}
\left [\, q^{i}, Y(u,z)\, \right] 
&=& 
\frac{1}{k!}\left [\, q^{i}, 
:\delz^{k}Y(a^{j},z)Y(v,z):\, \right]
\\
&=& 
\frac{1}{k!}\left [\, q^{i},\delz^{k}Y(a^{j},z)_{-}\, \right]Y(v,z)
\\
&&
+ \frac{1}{k!} :\delz^{k}Y(a^{j},z)Y(X^{i}_{+}(z)v,z):
\\
&=& 
Y\left(X^{i}_{+}(z)a^{j}(-k-1)v,z\right).
\end{eqnarray*} 
using $\left [\, q^{i},\delz^{k}Y(a^{j},z)_{+}\, \right]=0$ and
\begin{equation*}
\frac{1}{k!}\left [\, q^{i},\delz^{k}Y(a^{j},z)_{-}\, \right]=(-z)^{-k-1}\delta_{i,j}=\left[\, X^{i}_{+}(z),a^{j}(-k-1)\, \right].
\end{equation*} 
Hence \eqref{Yuconj2} holds by induction in $m$. The general result \eqref{Yuconj} follows on exponentiating and using $Y_{+}(\alpha,-z)=e^{-\alpha \cdot X_{+}(z)}$. $\square$
\medskip

We may now construct the creative intertwiner in much the same way as for a lattice vertex operator e.g.  \cite{K}:
\begin{theorem}
\label{theor_calY}
The creative intertwiner for $u\otimes e^{\alpha}\in M_{\alpha}$  for any $\alpha\in \mathbb{C}^{l}$ is given by
\begin{eqnarray}
{\cal Y}(u\otimes e^{\alpha},z)&=& e^{\alpha \cdot q}c_{\alpha} Y_{-}(\alpha, z)  
Y(u ,z)Y_{+}(\alpha,z)z^{\alpha(0)},
\label{Yalpha}
\end{eqnarray}
where $c_{\alpha}= \mathrm{Id}_{M^{l}}$.
\end{theorem}
\medskip
\textbf{Proof.} 
Using Proposition \ref{conjprop1} and skew-symmetry \eqref{skewW} we find that
\begin{eqnarray*}
{\cal Y}(u\otimes e^{\alpha},z)v &=& e^{zL(-1)}Y(v,-z)(u\otimes e^{\alpha})\\
&=& e^{zL(-1)}Y(v,-z)e^{\alpha \cdot q}u\\
&=& e^{\alpha \cdot q}Y_{-}(\alpha,z)e^{zL(-1)}Y(Y_{+}(\alpha,z)v,-z)u\\
&=& e^{\alpha \cdot q}Y_{-}(\alpha,z)Y(u,z)Y_{+}(\alpha,z)v,
\end{eqnarray*}
 for all $u,v\in M^{l}$. This implies that
\begin{equation}
{\cal Y}(u\otimes e^{\alpha},z)=e^{\alpha \cdot q} Y_{-}(\alpha, z)Y(u,z)  
Y_{+}(\alpha,z)b_{\alpha}(z),
\label{Ybal}
\end{equation} 
where $b_{\alpha}(z)\in \mathrm{End}(M^{l})[[z,z^{-1}]]$ with 
\begin{equation}
b_{0}(z)=\mathrm{Id}_{M^{l}},\quad [\,a^{i}(n),b_{\alpha}(z)\,]=0,\quad b_{\alpha}(z)v =v,
\label{balphav}
\end{equation}
for all $v\in M^{l}$. 
Translation \eqref{YL(-1)YInt} and $L(-1)(\vac\otimes e^{\alpha})= \alpha\cdot a\otimes e^{\alpha}$ imply
\begin{eqnarray*}
\delz {\cal Y}(\vac\otimes e^{\alpha},z)&=&{\cal Y}(\alpha\cdot a\otimes e^{\alpha},z).
\end{eqnarray*}
Using \eqref{Ybal} we find that $z\delz b_{\alpha}(z)=\alpha(0) b_{\alpha}(z)$ giving
\begin{equation}
b_{\alpha}(z)=c_{\alpha}z^{\alpha(0)},
\label{b_alpha}
\end{equation}
for $z$ independent operator $c_{\alpha}$. Applying \eqref{balphav} we conclude that $c_{\alpha}= \mathrm{Id}_{M^{l}}$ and hence the result follows. $\square$


\section{Generalized Vertex Operator Algebra for Heisenberg Intertwiners}
\label{GVOA}
\subsection{A Generalized Vertex Operator Algebra} \label{sub_GVOA}
Let ${\cal M}=\oplus_{\beta\in \mathbb{C}^l}M_{\beta}$, the direct sum of all the irreducible modules for $M^{l}=M_{0}$. Using \eqref{analpha}, the creative intertwining operator ${\cal Y}(u\otimes e^{\alpha},z)$ has a natural extension to an intertwiner vertex operator in $\mathrm{Lin}({\cal M},{\cal M})[[z,z^{-1}]]$ 
where now $c_{\alpha}$ acts on $M_{\beta}$ as a scalar
\begin{equation*}
c_{\alpha}u\otimes e^{\beta}=\epsilon(\alpha,\beta)u\otimes e^{\beta},
\end{equation*}
for $\epsilon(\alpha,\beta)\in \mathbb{C}^{\times}$ and $\epsilon(\alpha,0)=\epsilon(0,\alpha)=1$. 
As for lattice VOAs (e.g.  \cite{FLM}, \cite{K}), we define a cocycle system over $\mathbb{C}^{l}$ as an additive group. Define 
\begin{equation}
e^{\alpha}=e^{\alpha \cdot q} c_{\alpha},
\label{expal}
\end{equation}
for all $\alpha \in \mathbb{C}^{l}$ so that 
\begin{eqnarray}
{\cal Y}(u\otimes e^{\alpha},z)&=& e^{\alpha} Y_{-}(\alpha, z)  
Y(u ,z)Y_{+}(\alpha,z)z^{\alpha(0)}.
\label{Yalphaeal} 
\end{eqnarray}
We assume that the operators \eqref{expal} satisfy an associative algebra for 2-cocycle $\epsilon(\alpha,\beta)$ such that
\begin{equation}
e^{\alpha}e^{\beta}=\epsilon(\alpha,\beta)e^{\alpha+\beta},\quad e^{0}=1.
\label{epsab}
\end{equation}  
Associated to the cocycle system is the commutator function $C(\alpha,\beta)$ with
\begin{equation*}
e^{\alpha}e^{\beta}=C(\alpha,\beta)e^{\beta}e^{\alpha},
\end{equation*}
where
\begin{equation}
C(\alpha,\beta)=C(\beta,\alpha)^{-1}=\frac{\epsilon(\alpha,\beta)}{\epsilon(\beta,\alpha)}.
\label{Cab}
\end{equation}
Associativity implies $C(\alpha,\beta)$ is skewsymmetric and bilinear: 
\begin{eqnarray}
C(\alpha,\beta)&=&C(\beta,\alpha)^{-1},\notag\\
C(\alpha+\beta,\gamma)&=&C(\alpha,\gamma)C(\beta,\gamma),\notag\\
C(\alpha,\beta+\gamma)&=&C(\alpha,\beta)C(\alpha,\gamma).
\label{Cbil}
\end{eqnarray} 
\begin{example}
Suppose $l=2m$ and let $\alpha=(\alpha^{1},\alpha^{2})$ and 
$\beta=(\beta^{1},\beta^{2})$ for $\alpha^{i},\beta^{j}\in \mathbb{C}^{m}$. Then 
\begin{equation*}
C(\alpha,\beta)=\zeta^{\alpha^{1}\cdot\beta^{2}-\alpha^{2}\cdot\beta^{1}},
\end{equation*} 
for any $\zeta\in\mathbb{C}^{\times}$ satisfies \eqref{Cbil}.
\end{example}
\begin{lemma}
\label{Lem_coc}
The cocycle factors $\epsilon(\alpha,\beta)$ can be chosen such that $\epsilon(\alpha,-\alpha)=1$ for all $\alpha\in\mathbb{C}^{l}$.
\end{lemma}
\noindent \textbf{Proof.} Apply associativity to $e^{\alpha}e^{-\alpha}e^{\alpha}$ to find 
\begin{equation}
\epsilon(\alpha,-\alpha)=\epsilon(-\alpha,\alpha),\quad C(\alpha,-\alpha)=1.
\label{esym}
\end{equation}
As for lattice VOAs (e.g. \cite{K}), we may redefine $e^{\alpha}$ to be $\epsilon_{\alpha}e^{\alpha}$ for the same commutator function $C(\alpha,\beta)$ for any $\epsilon_{\alpha}\in \mathbb{C}^{\times}$ with $\epsilon_{0}=1$ where $\epsilon(\alpha,\beta)$ is redefined as  $\epsilon_{\alpha}\epsilon_{\beta}\epsilon_{\alpha+\beta}^{-1}\epsilon(\alpha,\beta)$.
Define an ordering on $\zeta\in\mathbb{C}$ with $\zeta>0$ if  $\Re(\zeta)>0$ or if 
$\Re(\zeta)=0$ and $\Im(\zeta)>0$. 
Choose 
\begin{equation}
\epsilon_{\alpha}=\left\{ 
\begin{array}{ll}
	\epsilon(\alpha,-\alpha)^{-1} & \mbox{ if } \alpha^{1}>0 \mbox{ or if } \alpha^{1}=\ldots = \alpha^{m-1}=0 \mbox{ and } \alpha^{m}>0, \\
	1 & \mbox{ otherwise. } 
\end{array}
\right.
\label{eq:}
\end{equation}
Hence $\epsilon(\alpha,-\alpha)$ is redefined as unity. $\square$   
\medskip

We also define the operator  \cite{Li1},  \cite{DLM}
\begin{equation}
\Delta(\alpha,z)=z^{\alpha(0)}Y_{+}(\alpha,-z).
\label{Delta}
\end{equation}
Using \eqref{Yuconj} and the cocycle structure \eqref{epsab} and \eqref{Cab} we immediately find: 
\begin{lemma} 
For all $\beta\in \mathbb{C}^{l}$ and $u\otimes e^{\alpha}\in M_{\alpha}$ 
\begin{equation}
(e^{\beta})^{-1}{\cal Y}(u\otimes e^{\alpha},z) e^{\beta}=C(\alpha,\beta){\cal Y}(\Delta(\beta,z)(u\otimes e^{\alpha}),z).  
\label{calYconj}
\end{equation}
\end{lemma}

The operators \eqref{Yalphaeal} with the above cocycle structure satisfy a natural extension from rational to complex parameters of the notion of a Generalized VOA as described in Chapt.~9 of  \cite{DL} and utilized in  \cite{DLM}. In this case, the operators \eqref{Yalphaeal} obey creativity and translation with Heisenberg Virasoro vector and satisfy a generalized Jacobi identity as follows
\begin{theorem}
\label{GenVOA}
The vertex operators ${\cal Y} (u\otimes e^{\alpha}, z )  \in\mathrm{Lin}({\cal M},{\cal M})[[z,z^{-1}]]$ satisfy the generalized Jacobi identity  
\begin{eqnarray}
\notag
& & z_0^{-1} \left( \frac{z_1 - z_2}{z_0}\right)^{ 
-\alpha \cdot\beta
}  
\delta\left( \frac{z_1 - z_2}{z_0}\right)  
\; {\cal Y} (u\otimes e^{\alpha}, z_1 ) \; 
{\cal Y} (v \otimes e^{\beta}, z_2)   
\\
\notag
& & \; 
  - C(\alpha,\beta)z_0^{-1}
\left( \frac{z_2 - z_1}{z_0}\right)^{ -\alpha\cdot \beta 
}  
 \delta\left( \frac{z_2 - z_1}{-z_0}\right)
 \;{\cal Y} (v \otimes e^{\beta}, z_2)  \; {\cal Y} (u\otimes e^{\alpha}, z_1 ) 
\\
\label{vogja}
&& 
\;
= z_2^{-1}  
 \;   
\delta\left( \frac{z_1 - z_0}{z_2}\right)
{\cal Y} ( {\cal Y} (u \otimes e^{\alpha}, z_0)
 (v\otimes e^{\beta}), z_2)
\left( \frac{z_1 - z_0}{z_2}\right)^{\alpha(0)},  
\end{eqnarray} 
for all $u\otimes e^{\alpha},v\otimes e^{\beta}\in{\cal M}$ with cocycle structure \eqref{epsab} and \eqref{Cab}.
\end{theorem}
\medskip
\noindent {\bf Proof.} 
The proof is similar to that of Theorem~3.5 of  \cite{DLM}. Using \eqref{Yuconj}
\begin{eqnarray*}
&&{\cal Y} (u\otimes e^{\alpha}, z_1 ){\cal Y} (v \otimes e^{\beta}, z_2) (w\otimes e^{\gamma})\\
&&=
z_{1}^{\alpha\cdot(\beta+\gamma)}z_{2}^{\beta\cdot\gamma}e^{\alpha}e^{\beta} Y_{-}(\alpha,z_1)Y(Y_{+}(\beta,-z_1)u,z_1)
\\
&&\quad .Y_{+}(\alpha,z_1)Y_{-}(\beta,z_2)Y(v,z_2)Y_{+}(\beta,z_2)
(w\otimes e^{\gamma})
\\
&&=
z_{1}^{\alpha\cdot\gamma}z_{2}^{\beta\cdot\gamma}
\left(z_{1}-z_{2} \right)^{\alpha \cdot\beta}
e^{\alpha}e^{\beta} 
Y_{-}(\alpha,z_1)Y(Y_{+}(\beta,-z_1)u,z_1)
\\
&&\quad .Y_{-}(\beta,z_2)Y_{+}(\alpha,z_1)Y(v,z_2)Y_{+}(\beta,z_2)
(w\otimes e^{\gamma})
\\
&&=
z_{1}^{\alpha\cdot\gamma}z_{2}^{\beta\cdot\gamma}\left(z_{1}-z_{2} \right)^{\alpha\cdot \beta}
e^{\alpha}e^{\beta}
Y_{-}(\alpha,z_1)Y_{-}(\beta,z_2)Y(Y_{+}(\beta,-z_1+z_{2})u,z_1)
\\
&&\quad .Y(Y_{+}(\alpha,z_1-z_{2})v,z_2)
Y_{+}(\alpha,z_1)Y_{+}(\beta,z_2)
(w\otimes e^{\gamma}),
\end{eqnarray*}
using \eqref{Ypmcom} and \eqref{Yuconjminus}.
Thus 
\begin{eqnarray*}
&&z_0^{-1} \left( \frac{z_1 - z_2}{z_0}\right)^{-\alpha\cdot \beta}  
\delta\left( \frac{z_1 - z_2}{z_0}\right)  
{\cal Y} (u\otimes e^{\alpha}, z_1 ) 
{\cal Y} (v \otimes e^{\beta}, z_2)
(w\otimes e^{\gamma})
\\
&&= 
z_{0}^{\alpha\cdot\beta}z_{1}^{\alpha\cdot\gamma}z_{2}^{\beta\cdot\gamma}
e^{\alpha}e^{\beta}Y_{-}(\alpha,z_1)Y_{-}(\beta,z_2)\\
&&
\quad
.z_0^{-1}  
\delta\left( \frac{z_1 - z_2}{z_0}\right)
Y(Y_{+}(\beta,-z_{0})u,z_1)Y(Y_{+}(\alpha,z_{0})v,z_2)
\\
&&
\quad
.Y_{+}(\alpha,z_1)Y_{+}(\beta,z_2)
(w\otimes e^{\gamma}),
\end{eqnarray*}
using $\delta\left( \frac{z_1 - z_2}{z_0}\right)(z_{1}-z_{2})^{n}=\delta\left( \frac{z_1 - z_2}{z_0}\right)z_{0}^n$ for any integer $n$. 
Similarly, we find 
\begin{eqnarray*}
&&C(\alpha,\beta)z_0^{-1} \left( \frac{z_2 - z_1}{z_0}\right)^{-\alpha \cdot\beta}  
\delta\left( \frac{z_2 - z_1}{-z_0}\right)  
{\cal Y} (v \otimes e^{\beta}, z_2)
{\cal Y} (u\otimes e^{\alpha}, z_1 )
(w\otimes e^{\gamma})
\\
&&= 
z_{0}^{\alpha\cdot\beta}z_{1}^{\alpha\cdot\gamma}z_{2}^{\beta\cdot\gamma}
e^{\alpha}e^{\beta}
Y_{-}(\alpha,z_1)Y_{-}(\beta,z_2) \\
&&\quad
.z_0^{-1}\delta\left( \frac{z_2 - z_1}{-z_0}\right)
Y(Y_{+}(\alpha,z_{0})v,z_2)Y(Y_{+}(\beta,-z_{0})u,z_1)
\\
&&\quad .Y_{+}(\alpha,z_1)Y_{+}(\beta,z_2)
(w\otimes e^{\gamma}).
\end{eqnarray*}
Hence on applying the Jacobi identity for the Heisenberg VOA, the left hand side of \eqref{vogja} applied to $w\otimes e^{\gamma} $  gives
\begin{eqnarray}
&&z_2^{-1} 
\delta\left( \frac{z_1 - z_0}{z_2}\right) z_{0}^{\alpha\cdot\beta}z_{1}^{\alpha\cdot\gamma}z_{2}^{\beta\cdot\gamma}
e^{\alpha}e^{\beta}
Y_{-}(\alpha,z_1)Y_{-}(\beta,z_2)\notag
\\
&&.Y(B,z_2)Y_{+}(\alpha,z_1)Y_{+}(\beta,z_2)
(w\otimes e^{\gamma}).
\label{LHS}
\end{eqnarray}
for $B=Y(Y_{+}(\beta,-z_{0})u, z_{0} )Y_{+}(\alpha,z_{0})v$.
\medskip

In a similar way 
\begin{eqnarray*}
&&
{\cal Y} \left(
{\cal Y} (u\otimes e^{\alpha}, z_0 )(v \otimes e^{\beta}), z_2
\right) 
(w\otimes e^{\gamma})
\\
&&=
z_{0}^{\alpha\cdot\beta}\epsilon(\alpha,\beta)
{\cal Y}\left(
Y_{-}(\alpha,z_{0})
B\otimes e^{\alpha+\beta}, z_{2} 
\right) 
(w\otimes e^{\gamma})
\\
&&=
z_{0}^{\alpha\cdot\beta}z_{2}^{(\alpha+\beta)\cdot\gamma}e^{\alpha}e^{\beta}
Y_{-}(\beta,z_{2})Y_{-}(\alpha,z_{2})
Y\left( 
Y_{-}(\alpha,z_{0})B
, z_{2} 
\right)
\\
&&\quad
.Y_{+}(\alpha,z_{2})Y_{+}(\beta,z_{2}) 
(w\otimes e^{\gamma}).
\end{eqnarray*}
Employing \eqref{YYconj} the right hand side of \eqref{vogja} applied to $w\otimes e^{\gamma} $  therefore gives
\begin{eqnarray}
&&z_2^{-1}  \left( \frac{z_1 - z_0}{z_2}\right)^{\alpha \cdot\gamma}  
\delta\left(\frac{z_1 - z_0}{z_2}\right)
z_{0}^{\alpha\cdot\beta}\left(z_{2}+z_{0}\right)^{\alpha\cdot\gamma}z_{2}^{\beta\cdot\gamma}
e^{\alpha}e^{\beta}
\notag
\\
&& .Y_{-}(\beta,z_{2})Y_{-}(\alpha,z_{0}+z_{2})
Y(B, z_{2})Y_{+}(\alpha,z_{2}+z_{0})
Y_{+}(\beta,z_{2}) 
(w\otimes e^{\gamma}).
\label{RHS}
\end{eqnarray}
Finally, using the identity 
\begin{equation*}
\delta\left( \frac{z_{1}-z_{0}}{z_{2}} \right)
\left( \frac{z_{1}-z_{0}}{z_{2}}\right) ^{\kappa}
=
\delta\left(  \frac{z_{1}-z_{0}}{z_{2}}  \right)
\left( \frac{z_{2}+z_{0}}{z_{1}}\right) ^{-\kappa},
\end{equation*}
for $\kappa\in \mathbb{C}$,
we find that \eqref{LHS} and \eqref{RHS} are equal. Thus the theorem holds.
$\square$

\subsection{Skew-Symmetry and an Invariant Form} \label{sub_skew}
In order to formulate a generalization of skew-symmetry \eqref{skew} and \eqref{skewW} applicable to ${\cal Y} (u\otimes e^{\alpha}, z )$ we firstly define for formal parameters $z, \kappa$  
\begin{equation}
(-z)^{\kappa}=e^{i\pi N\kappa }z^{\kappa},
\label{minz}
\end{equation}
where $N$ is an odd integer parameterizing the formal branch cut.
\begin{lemma}
The operators ${\cal Y} (u\otimes e^{\alpha}, z)$ satisfy the skew-symmetry property
\begin{equation}
\label{skew_symm}
{\cal Y}(u\otimes e^{\alpha}, z) (v\otimes e^{\beta})= 
e^{-i\pi N \alpha \cdot\beta } 
C(\alpha,\beta)  
e^{zL(-1)} {\cal Y}(v\otimes e^{\beta}, -z)  (u\otimes e^{\alpha}).  
\end{equation}
\end{lemma} 
\noindent
\textbf{Proof.}
Using   \eqref{skewW}, \eqref{L1conj} and \eqref{calYconj} we have ${\cal Y}(u\otimes e^{\alpha}, z) (v\otimes e^{\beta})$ is given by 
\begin{eqnarray*}
\label{skew_symm_2}
&&
{\cal Y}(u\otimes e^{\alpha}, z)  e^{\beta} v
\\
&& =C(\alpha,\beta) e^{\beta }  {\cal Y}( \Delta(\beta,z) (u \otimes e^{\alpha}), z) v
\\
&& = C(\alpha,\beta) e^{\beta}  
 e^{z L(-1)} Y ( v , -z) \Delta(\beta,z)(u \otimes e^{\alpha})  
 \\
 && = C(\alpha,\beta) e^{z L(-1)} e^{\beta}  Y_{-}(\beta,-z)
 Y ( v , -z) Y_{+}(\beta,-z)z^{\alpha\cdot\beta} (u \otimes e^{\alpha})  
 \\
 && = e^{-i\pi N \alpha \cdot\beta }C(\alpha,\beta)
e^{z L(-1)} {\cal Y} ( v\otimes e^{\beta} , -z) (u \otimes e^{\alpha}).\  \square 
\end{eqnarray*}

\medskip
We next introduce an invariant form $\langle \, ,\rangle$ on ${\cal M}$  
associated with the M\"{o}bius map  \cite{FHL}, \cite{Li2}, \cite{S}, \cite{TZ}
\begin{equation}
\left(
\begin{array}{cc}
0 & \lambda\\
e^{i\pi N }\lambda^{-1} & 0\\	
\end{array}
\right)
:z\mapsto \frac{\lambda^{2}}{e^{i\pi N }z},
 \label{eq: gam_lam}
\end{equation}
for $\lambda\neq 0$ and with $N$ of \eqref{minz}.
Usually for a VOA we have $e^{-i\pi N }=-1$ and one takes $\lambda=\pm \sqrt{-1}$ since only integral powers of formal parameters occur. In the present case we define the adjoint of the vertex operator
${\cal Y}(u\otimes e^{\alpha},z) $ by
\begin{eqnarray}
{\cal Y}^\dagger\left(u\otimes e^{\alpha},z\right) &=&
{\cal Y}\left(
e^{-z\lambda^{-2}L(1)} 
\left(\frac{\lambda}{e^{i\pi N }z}\right)^{2 L(0)}
(u\otimes e^{\alpha})
,\frac{\lambda^{2}}{e^{i\pi N }z}\right). 
\label{eq: calY_dag}
\end{eqnarray}   
In particular, for a Heisenberg generating vector $a^{i}$ we have
\begin{equation}
Y(a^{i},z)^{\dagger }=\sum_{n}{a^{i}}^{\dagger }(n)z^{-n-1},
\quad {a^{i}}^{\dagger }(n)=(-1)^{n+1}\lambda^{2n} a^{i}(-n).   
\label{adagger}
\end{equation}
This implies 
\begin{equation}
\label{ypmdagger}
Y_\pm^{\dagger} (\alpha, z)= Y_\mp \left(\alpha, -\frac{\lambda^2}{z}\right).    
\end{equation}
We also note that $e^{-z\lambda^{-2}L(1)} 
\left(\frac{\lambda}{e^{i\pi N } z}\right)^{2 L(0)}
(u\otimes e^{\alpha})$ is given by
\begin{eqnarray*}
&&\left(e^{-i\pi N }\frac{\lambda}{z}\right)^{\alpha^2} 
e^{-z\lambda^{-2}L(1)} 
\left(\left(-\frac{\lambda}{z}\right)^{2 L(0)}u\right)\otimes e^{\alpha}
\\
&&\quad =\left(e^{-i\pi N }\frac{\lambda}{z}\right)^{\alpha^2} 
e^{-z\lambda^{-2}L(1)}  e^{\alpha \cdot q} 
\left(-\frac{\lambda}{z}\right)^{2 L(0)}u
\\
&&\quad = \left(e^{-i\pi N }\frac{\lambda}{z}\right)^{\alpha^2} e^{\alpha \cdot q} 
Y_{+}\left(\alpha,\frac{\lambda^{2}}{z}\right)
e^{-z\lambda^{-2}L(1)}\left(-\frac{\lambda}{z}\right)^{2 L(0)}u
\\
&& \quad=\left(e^{-i\pi N }\frac{\lambda}{ z}\right)^{\alpha^2}
\left(Y_{+}\left(\alpha,\frac{\lambda^{2}}{z}\right)
e^{-z\lambda^{-2}L(1)}
\left(-\frac{\lambda}{z}\right)^{2 L(0)}u\right)
\otimes e^{\alpha},
\end{eqnarray*} 
where $\alpha^2=\alpha\cdot \alpha$ and using \eqref{Lplus1conj}. 
Hence we find
\begin{eqnarray}
&&{\cal Y}^\dagger\left(u\otimes e^{\alpha},z\right)\, (w\otimes e^{\gamma})
\notag
\\
&=&\left(e^{-i\pi N }\frac{\lambda}{z}\right)^{\alpha^2}
\left(e^{-i\pi N }\frac{\lambda^2}{ z}\right)^{\alpha\cdot\gamma} 
e^{\alpha } Y_{-}\left(\alpha,-\frac{\lambda^{2}}{z}\right)
\notag
\\
&&
.Y\left(
Y_{+}\left(\alpha,\frac{\lambda^{2}}{z}\right)
e^{-z\lambda^{-2}L(1)}
\left(-\frac{\lambda}{z}\right)^{2 L(0)}u,
-\frac{\lambda^{2}}{z}
\right)\notag
\\
&&
.
Y_{+}\left(\alpha,-\frac{\lambda^{2}}{z}\right)
 (w\otimes e^{\gamma}) 
\notag
\\
&=&
{\lambda}^{2\alpha\cdot(\alpha+\gamma)-\alpha^2}(e^{i\pi N} z)^{-\alpha\cdot(\alpha+\gamma)}
Y_{-}\left(\alpha,-\frac{\lambda^{2}}{z}\right)
\notag
\\
&&
.Y\left(
e^{-z\lambda^{-2}L(1)}
\left(-\frac{\lambda}{z}\right)^{2 L(0)}u,
-\frac{\lambda^{2}}{z}
\right)\notag
\\
&&
.
Y_{+}\left(\alpha,-\frac{\lambda^{2}}{z}\right)
e^{\alpha }  (w\otimes e^{\gamma}) 
\notag
\\
&=&
z^{-\alpha(0)}
Y_{-}\left(\alpha,-\frac{\lambda^{2}}{z}\right)
Y^{\dagger}(u,z)Y_{+}\left(\alpha,-\frac{\lambda^{2}}{z}\right)
\notag
\\
&&.
e^{-i\pi N \alpha(0)}
\lambda^{2 \alpha(0)-\alpha^2}
e^{\alpha }  (w\otimes e^{\gamma}).
\notag
\end{eqnarray}
Thus using \eqref{adagger} and \eqref{ypmdagger} we have
\begin{equation}
{\cal Y}^\dagger\left(u\otimes e^{\alpha},z\right)
=z^{\alpha(0)^{\dagger}}
Y_{+}^{\dagger}(\alpha,z)
Y^{\dagger}(u,z)Y_{-}^{\dagger}(\alpha,z)
{e^{\alpha }}^{\dagger},
 \label{Ydag2}
\end{equation}
where we define 
\begin{equation}
{e^{\alpha }}^{\dagger}= 
e^{-i\pi N \alpha(0)}
\lambda^{2 \alpha(0)-\alpha^2}e^{\alpha }. 
\label{eq:qdag}
\end{equation}

\medskip

\begin{definition} A bilinear form $\langle \, ,\rangle$ on ${\cal M}$ is said to be   
invariant if for all $u\otimes e^{\alpha}$, $v\otimes e^{\beta} $, $w\otimes e^{\gamma} \in {\cal M}$ we have
\begin{equation}
\langle\, {\cal Y}(u\otimes e^{\alpha},z) (v\otimes e^{\beta})
, w\otimes e^{\gamma}\, \rangle 
= 
e^{-i\pi N \alpha\cdot\beta}
 C(\alpha,\beta)
\langle\, v\otimes e^{\beta}, \; {\cal Y}^{\dagger}(u\otimes e^{\alpha},z)( w\otimes e^{\gamma})\, \rangle. 
\label{eq: inv bil form}
\end{equation}
\end{definition} 
We choose the normalization $\langle\, {\bf 1},{\bf 1} \, \rangle=1$. 
For $\alpha=\beta=\gamma=0$ \eqref{eq: inv bil form} reverts to the standard 
definition\footnote{up to an additional $\lambda$ dependence arising from definition  for the adjoint in \eqref{adagger}.} of an invariant form on the Heisenberg VOA $M^{l}$ which is unique, symmetric and invertible  \cite{Li2}. In general, we have
\begin{proposition}
\label{forma_prop}
The bilinear form $\langle \, ,\rangle$ on $\cal{M}$ is unique, symmetric and invertible with 
\begin{equation}
\langle\, v\otimes e^{\beta}, w\otimes e^{\gamma} \, \rangle=
\epsilon(\beta,-\beta)
\lambda^{-\beta^2}
\delta_{\beta+\gamma,0}
\langle\,  v, w\, \rangle.
\label{albet}
\end{equation}
\end{proposition}
\begin{remark}
From Lemma~\ref{Lem_coc} we note that we may choose $\epsilon(\beta,-\beta)=1$. 
\end{remark}
\noindent{\bf Proof.} 
Since $\langle\, a(0)(v\otimes e^{\beta}), w\otimes e^{\gamma}\, \rangle=-\langle\, v\otimes e^{\beta}, a(0)(w\otimes e^{\gamma})\, \rangle$ it follows that 
\begin{eqnarray*}
\langle\, v\otimes e^{\beta}, w\otimes e^{\gamma} \, \rangle=
\delta_{\beta+\gamma,0}\langle\, v\otimes e^{\beta}, w\otimes e^{-\beta}\, \rangle.
\end{eqnarray*} 
Applying \eqref{eq:qdag} we obtain
\begin{eqnarray*}
\langle\, v\otimes e^{\beta}, w\otimes e^{-\beta}\, \rangle
&=& 
\langle\, e^{\beta}v, w\otimes e^{-\beta}\, \rangle\\
&=&  
\langle\, v,  {e^{\beta}}^{\dagger}(w\otimes e^{-\beta})\, \rangle\\
&=&  
\epsilon(\beta,-\beta)\lambda^{-\beta^2}\langle\, v, w\, \rangle.  
\end{eqnarray*}
But $\epsilon(\beta,-\beta)=\epsilon(-\beta,\beta)$ from \eqref{esym} and $\langle\, v, w\, \rangle$ is symmetric, unique and invertible so the result holds. $\square$  

\begin{remark}
The necessity for the external factors on the right hand side of \eqref{eq: inv bil form} is apparent when we consider 
\begin{eqnarray*}
\langle\, e^{\alpha} (v\otimes e^{\beta})
, w\otimes e^{\gamma}\, \rangle 
&=& 
\epsilon(\alpha,\beta) \langle\, v\otimes e^{\alpha+\beta}
, w\otimes e^{\gamma}\, \rangle\\
&=&  \epsilon(\alpha,\beta) \epsilon(\alpha+\beta,-\alpha-\beta)\lambda^{-(\alpha+\beta)^2}
\delta_{\alpha+\beta+\gamma,0}
\langle\, v, w\, \rangle.
\end{eqnarray*}
On the other hand \eqref{eq: inv bil form} implies that this is also given by
\begin{eqnarray*}
&& 
e^{-i\pi N \alpha\cdot\beta}C(\alpha,\beta)
\langle\, v\otimes e^{\beta}
, {e^{\alpha}}^{\dagger} (w\otimes e^{\gamma})\, \rangle 
\\
&&= 
e^{-i\pi N \alpha\cdot\beta}C(\alpha,\beta)
e^{-i\pi N \alpha\cdot(\alpha+\gamma)}
\lambda^{\alpha^2+2\alpha\cdot\gamma}
\epsilon(\alpha,\gamma)
\langle\, v\otimes e^{\beta}
, w\otimes e^{\alpha+\gamma}\, \rangle 
\\
&&=  C(\alpha,\beta)\epsilon(\alpha,-\alpha-\beta)\epsilon(\beta,-\beta) 
\lambda^{-(\alpha+\beta)^2}
\delta_{\alpha+\beta+\gamma,0}
\langle\, v, w\, \rangle.
\end{eqnarray*}
The equality of these expressions is equivalent to the identity 
$(e^{\alpha}e^{\beta})e^{-\alpha -\beta}=C(\alpha,\beta)e^{\beta}(e^{\alpha}e^{-\alpha -\beta})$.
\end{remark}

\section{Lattice Vertex Operator Superalgebras}
\label{LVOSA}
In this section we apply Theorem~\ref{GenVOA} to the example of a lattice VOSA
$V_{L}$ for an integral
Euclidean rank $l$ lattice $L$. We construct the $g$-twisted module for a $V_{L}$
automorphism $g$ generated by a Heisenberg vector in terms of Heisenberg modules so that
the twisted Jacobi identity \eqref{gJac} is satisfied as a consequence of Theorem~\ref{GenVOA}. 
The relationship between this and Li's construction \cite{Li1} for a  $g$-twisted module is discussed.   
We also consider a generalization from rational to complex parameterized twisted $V_{L}$ modules VOSA of a related generalized VOA discussed in \cite{DLM}.

Let $L$ be a Euclidean lattice of rank $l$ and define $V_{L}=\oplus_{\mu\in L}M_{\mu}$ with standard cocycle commutator e.g. \cite{K}
\begin{equation}
C(\mu_{1},\mu_{2})=(-1)^{\mu_{1}\cdot\mu_{2}+\mu_{1}^2\mu_{2}^2},
\label{Clattice}
\end{equation}  
for $\mu_{1},\mu_{2}\in L$.
Define parity on $V_{L}$ by $p(u\otimes e^{\mu})=\mu^{2}\mod 2$ for $u\otimes e^{\mu}\in V_{L}$. Then
Theorem~\ref{GenVOA} implies that $(V_{L},{\cal Y},\vac,\omega)$ is a VOA for $L$ even and a VOSA for $L$ odd with invertible invariant form \eqref{eq: inv bil form} obeying
\begin{equation*}
\langle\, {\cal Y}(u\otimes e^{\mu_{1}},z) (v\otimes e^{\mu_{2}})
, w\otimes e^{\mu_{3}}\, \rangle 
= 
(-1)^{\mu_{1}^2\mu_{2}^2}
\langle\, v\otimes e^{\mu_{2}}, {\cal Y}^{\dagger}(u\otimes e^{\mu_{1}}),z)( w\otimes e^\mu_{3})\, \rangle. 
\end{equation*}

Consider the automorphism 
\begin{equation}
g=e^{-2\pi i \alpha(0)},
\label{galpha}
\end{equation}
generated by the Heisenberg vector $-\alpha\cdot a\in V_{L}$ for any $\alpha\in \mathbb{C}^{l}$. Clearly $M_{\mu}$ is a $g$ eigenspace for eigenvalue 
$e^{-2\pi i \alpha\cdot\mu}$. Let $V_{L+\alpha}=\oplus_{\mu\in L}M_{\mu+\alpha}$ so that 
$e^{\alpha}:V_{L}\rightarrow V_{L+\alpha}$. We then find using Definition~\ref{Vgmodule} that:
\begin{proposition}
\label{VLg}
$(V_{L+\alpha},{\cal Y})$ is a $g$-twisted $V_{L}$-module.
\end{proposition}
\textbf{Proof.} 
Theorem~\ref{GenVOA} with commutator \eqref{Clattice} implies that
\begin{eqnarray}
\notag
& & z_0^{-1}   
\delta\left( \frac{z_1 - z_2}{z_0}\right)  
{\cal Y} (u\otimes e^{\mu_{1}}, z_1 )
{\cal Y} (v \otimes e^{\mu_{2}}, z_2)
(w\otimes e^{\mu_{3}+\alpha})   
\\
\notag
& &  
 - (-1)^{\mu_{1}^2\mu_{2}^2}
 z_0^{-1}
 \delta\left( \frac{z_2 - z_1}{-z_0}\right)
{\cal Y} (v \otimes e^{\mu_{2}}, z_2) 
 {\cal Y} (u\otimes e^{\mu_{1}}, z_1 ) 
 (w\otimes e^{\mu_{3}+\alpha})
\\
\notag
&& 
= z_2^{-1}\left( \frac{z_1 - z_0}{z_2}\right)^{\mu_{1}\cdot \alpha}  
\delta\left( \frac{z_1 - z_0}{z_2}\right)
{\cal Y} ( {\cal Y} (u \otimes e^{\mu_{1}}, z_0)
 (v\otimes e^{\mu_{2}}), z_2)(w\otimes e^{\mu_{3}+\alpha}),
 \\
 &&
\label{YgYcal}  
\end{eqnarray} 
for $u\otimes e^{\mu_{1}}, v\otimes e^{\mu_{2}}\in V_{L}$ and $w\otimes e^{\mu_{3}+\alpha}\in V_{L+\alpha}$. The result holds on comparison with \eqref{gJac} where $\rho=-\mu_{1}\cdot \alpha$. $\square$

Using Lemma~\ref{calYconj} we immediately find that Proposition~\ref{VLg} is equivalent to Li's construction for the $g$-twisted module for $V_{L}$ \cite{Li1}:
\begin{corollary}
$(V_{L+\alpha},{\cal Y})\cong (V_{L},Y_{g})$ as $g$-twisted $V_{L}$-modules
where 
\begin{equation*}
Y_{g}(u\otimes e^{\mu},z)={\cal Y}(\Delta(\alpha,z)(u\otimes e^{\mu}),z).
\end{equation*}
\end{corollary}
For the Heisenberg basis and Virasoro vector we obtain $g$-twisted modes
\begin{eqnarray}
a_{g}^{i}(n)&=&(e^{\alpha})^{-1}a^{i}(n)e^{\alpha}
\notag
\\
&=&
a^{i}(n)+\alpha^{i}\delta_{n,0},
\label{agn}\\
\quad L_{g}(n)&=&(e^{\alpha})^{-1}L(n)e^{\alpha}
\notag
\\
&=&
\frac{1}{2}\sum_{m\in \mathbb{Z}}:a_{g}(n+m)\cdot a_{g}(-m):
\notag
\\
&=& 
L(n)+\alpha(n)+\frac{1}{2}\alpha^2\delta_{n,0},
\label{Lgn}
\end{eqnarray}
satisfying the original Heisenberg and central charge $l$ Virasoro algebras.
In particular, the grading is determined by $L_{g}(0)=L(0)+\alpha(0)+\frac{1}{2}\alpha^2$.

For each Heisenberg vector $\alpha\cdot a$ we may also construct a $\mathbb{C}$-graded VOSA $(V_{L},{\cal Y},\mathbf{1},\omega_{\alpha})$ with the original vector space and vertex operators but with a new \lq shifted\rq\ conformal vector $\omega_{\alpha}=\omega -\alpha(-2)\mathbf{1}$ \cite{MN}, \cite{DM}. The $\mathbb{C}$-grading is determined by $L_{\alpha}(0)=L(0)+\alpha(0)$ with twisted Virasoro modes 
\begin{equation*}
L_{\alpha}(n)=L(n)+(n+1)\alpha(n),
\end{equation*} 
satisfying the Virasoro algebra with central charge 
 $c_{\alpha}=l-12\alpha^2$. Thus we find  \cite{DM}
\begin{equation*}
L_{\alpha}(0)-\frac{c_{\alpha}}{24}=L_{g}(0)-\frac{l}{24},
\end{equation*}   
a fact that has been usefully exploited to relate shifted and twisted partition and $n$-point functions \cite{DM}, \cite{MTZ}.

\medskip

We conclude with a generalization of one of the main results of ref.~\cite{DLM} where a 
generalized VOA is constructed from an extension of a simple VOA $V$ by 
Heisenberg generated $g$-twisted modules with rational parameters via Li's construction.
Applying this approach to $V_{L}$ the 
generalized vertex operators act on $g$-twisted modules $\{V_{L+\alpha}\}$ in our notation. 
Applying Definition~3.3 of ref.~\cite{DLM} to the VOSA $V_{L}$ with $\psi_{\alpha}=e^{-\alpha\cdot q}$ we define the vertex operator
\begin{eqnarray}
Y_{\alpha}(u\otimes e^{\mu_{1}+\alpha},z)(v\otimes e^{\mu_{2}+\beta})
&=&
\psi_{-\alpha-\beta}Y_{-}(\alpha,z)
{\cal Y}(\psi_{\alpha}\Delta(\beta,z)(u\otimes e^{\mu_{1}+\alpha_{1}}),z)
\notag\\
&&.\Delta(\alpha,-z)\psi_{\beta}(v\otimes e^{\mu_{2}+\beta}),
\label{YalphaDLM}
\end{eqnarray}
for $\mu_{1},\mu_{2}\in L$. Since $\Delta(\alpha,-z)=Y_{+}(\alpha,z)(-z)^{\alpha(0)}$ 
we must employ the formal branch parameterization of \eqref{minz} to find
\begin{eqnarray*}
Y_{\alpha}(u\otimes e^{\alpha + \mu_{1}},z) (v\otimes e^{\mu_{2}+\beta})
&=&
e^{i\pi N \alpha \cdot \mu_{2}}
\epsilon(\mu_{1},\mu_{2})Y_{-}(\mu_{1}+\alpha,z)
Y(u,z)
\\
&&.Y_{+}(\mu_{1}+\alpha,z)z^{(\mu_{1}+\alpha)a(0)}(v\otimes e^{\mu_{2}+\beta})
\\
&=&
\frac{e^{i\pi N \alpha \cdot \mu_{2}}
\epsilon(\mu_{1},\mu_{2})}{\epsilon(\mu_{1}+\alpha,\mu_{2}+\beta)}
{\cal Y}(u\otimes e^{\mu_{1}+\alpha},z) (v\otimes e^{\mu_{2}+\beta}).
\end{eqnarray*}
Substituting into Theorem~\ref{GenVOA} results in
\begin{proposition}
\label{DLMGenVOA}
The vertex operators $Y_{\alpha}(u\otimes e^{\alpha + \mu},z) $ satisfy the generalized Jacobi identity  
\begin{eqnarray}
\notag
&  z_0^{-1} \left( \frac{z_1 - z_2}{z_0}\right)^
{\eta_{12}}
\delta\left( \frac{z_1 - z_2}{z_0}\right)  
Y_{\alpha_{1}}(u\otimes e^{\mu_{1}+\alpha_{1}}, z_1 )
Y_{\alpha_{2}} (v \otimes e^{\mu_{2}+\alpha_{2}}, z_2) (w\otimes e^{\mu_{3}+\alpha_{3}})&  
\\
\notag
& -C_{12}
  z_0^{-1}
\left( \frac{z_2 - z_1}{z_0}\right)^
{\eta_{12}}
 \delta\left( \frac{z_2 - z_1}{-z_0}\right)
Y_{\alpha_{2}}(v \otimes e^{\mu_{2}+\alpha_{2}}, z_2)  
Y_{\alpha_{1}} (u\otimes e^{\mu_{1}+\alpha_{1}}, z_1 ) 
(w\otimes e^{\mu_{3}+\alpha_{3}})& 
\\
\notag
&= z_2^{-1} 
\left( \frac{z_1 - z_0}{z_2}\right)^{-\eta_{13}}    
\delta\left( \frac{z_1 - z_0}{z_2}\right)
Y_{\alpha_{1}+\alpha_{2}}( Y_{\alpha_{1}} (u \otimes e^{\mu_{1}+\alpha_{1}}, z_0)
 (v\otimes e^{\mu_{2}+\alpha_{2}}), z_2) 
(w\otimes e^{\mu_{3}+\alpha_{3}})
&\\
&&\label{DLMvogja}\end{eqnarray} 
where 
\begin{eqnarray}
\eta_{12}&=&-\alpha_{1}\cdot \alpha_{2} -\mu_{1}\cdot \alpha_{2}-\mu_{2}\cdot \alpha_{1},
\label{eta12}
\\
\eta_{13}&=&-\alpha_{1}\cdot \alpha_{3} -\mu_{1}\cdot \alpha_{3}-\mu_{3}\cdot \alpha_{1},
\label{eta13}
\end{eqnarray}
and with commutator 
\begin{equation}
C_{12}=e^{i\pi N (\alpha_{1} \cdot \mu_{2}-\alpha_{2} \cdot \mu_{1})}(-1)^{\mu_{1}^2\mu_{2}^2}.
\label{C12}
\end{equation}
\end{proposition}
Choosing the branch\footnote{This branch choice is an unstated assumption in eqn.~(3.30)
 in the proof of Theorem~3.5 of ref.~\cite{DLM}} \eqref{minz} with $N=1$  Theorem~\ref{DLMGenVOA} extends Theorem~3.5
 of ref.~\cite{DLM} from rational to complex parametrized twisted modules of a VOSA $V_{L}$.
 Furthermore, it is clear that the commutator factor of $e^{i\pi N (\alpha_{1} \cdot \mu_{2}-\alpha_{2} \cdot \mu_{1})}$ arises solely from the branch choice made in \eqref{YalphaDLM}. In fact, we may modify  Defn.~3.3 of \cite{DLM} and \eqref{YalphaDLM} by replacing the $\Delta(\alpha,-z)$ operator by $Y_{+}(\alpha,z)z^{\alpha(0)}$ in order to define a new vertex operator  
\begin{eqnarray}
\widehat{Y}_{\alpha}(u\otimes e^{\mu_{1}+\alpha},z)(v\otimes e^{\mu_{2}+\beta})
&=&
\psi_{-\alpha-\beta}Y_{-}(\alpha,z)
{\cal Y}(\psi_{\alpha}\Delta(\beta,z)(u\otimes e^{\mu_{1}+\alpha_{1}}),z)
\notag\\
&&.Y_{+}(\alpha,z)z^{\alpha(0)}\psi_{\beta}(v\otimes e^{\mu_{2}+\beta}).
\label{Yalphahat}
\end{eqnarray}
These operators satisfy the generalized Jacobi identity \eqref{DLMvogja} with the standard lattice parity commutator $C_{12}=(-1)^{\mu_{1}^2\mu_{2}^2}$. 

\begin{thebibliography}{DLinM}

\bibitem{B} R. Borcherds, 
Vertex algebras, Kac-Moody algebras and the Monster, 
Proc.Natl.Acad.Sci. U.S.A. \textbf{83} (1986) 3068--3071.

\bibitem{FLM} I. Frenkel,   J. Lepowsky, A. Meurman,  
Vertex Operator Algebras and the Monster,
Academic Press, New York, 1988.

\bibitem{DLM} C. Dong, H. Li, G. Mason, 
Simple currents and extensions of vertex operator algebras,
Comm.Math.Phys. \textbf{180} (1996) 671--707.

\bibitem{Li1} H. Li, 
Local systems of twisted vertex operators, vertex operator superalgebras and twisted modules,
Contemp.Math. \textbf{193} (1996) 203--236.

\bibitem{DL} C. Dong,  J. Lepowsky,  
Generalized Vertex Algebras and Relative Vertex Operators,
 Progress in Mathematics, \textbf{112}, Birkhauser, Boston,  1993.

\bibitem{FHL} I. Frenkel,  Y-Z. Huang,  J. Lepowsky, 
On axiomatic approaches to vertex operator algebras and modules, 
Mem.Amer.Math.Soc. \textbf{104} no. 494 (1993).

\bibitem{K} V. Kac, 
Vertex Operator Algebras for Beginners,
University Lecture Series, Vol. 10, AMS, Boston, 1998.

\bibitem{LL}  J. Lepowsky,  H. Li,
Introduction to Vertex Operator Algebras and their Representations, 
Progr.Math. \textbf{227}, Birkhauser, Boston, 2004.

\bibitem{MT}  G. Mason, M.P. Tuite, 
Vertex operators and modular forms, 
\textit{A Window into Zeta and Modular Physics} eds. K.~Kirsten and F.~Williams, 
MSRI Publications \textbf{57} 183--278, Cambridge University Press, Cambridge, 2010.

\bibitem{X}  X. Xu, 
Introduction to Vertex Operator Superalgebras and their Modules,
Kluwer Academic, 1998. 

\bibitem{D} C. Dong,  
Twisted modules for vertex operator algebras associated with even lattices,
J.Alg. \textbf{165} (1993) 91--112.

\bibitem{DLinM} C. Dong, Z. Lin, G. Mason, 
On vertex operator algebras as $sl_{2}$-modules,
Arasu, K. T. (ed.) et al., 
\textit{Groups, Difference Sets, and the Monster},
Proceedings of a special research quarter, Columbus, OH, USA, Spring 1993,
Walter de Gruyter, Berlin, 1996,  
Ohio~State~Univ.Math.Res.Inst.Publ. \textbf{4} (1996) 349--362.  

\bibitem{Li2} H. Li, 
Symmetric invariant bilinear forms on vertex operator algebras, 
J.Pure.Appl.Alg. \textbf{96} (1994) 2790--297.

\bibitem{S} N. Scheithauer,  
Vertex algebras, Lie algebras and superstrings, 
J.Alg. \textbf{200} (1998) 363--403.

\bibitem{TZ} M.P. Tuite,  A. Zuevsky, 
Genus two partition and correlation functions for fermionic vertex operator superalgebras~I, Commun.Math.Phys. \textbf{306} (2011) 419--447. 

\bibitem{MN} A. Matsuo, K. Nagatomo,
A note on free bosonic vertex algebras and its conformal vector, 
J.Alg. \textbf{212} (1999) 395--418. 

\bibitem{DM} C. Dong,  G. Mason,  
Shifted vertex operator algebras,
Math.Proc.Camb.Philos.Soc. \textbf{141} (2006) 67--80.
 
\bibitem{MTZ} G. Mason, M.P Tuite, A. Zuevsky,
Torus n-point functions for $\mathbb{C}$-graded vertex operator superalgebras and continuous
fermion orbifolds,  
Commun.Math.Phys. \textbf{283} (2008) 305--342.  

\end{thebibliography}
\end{document}